\documentstyle[11pt,twoside]{article}
\include{graphicx}
\title{A note on the asymptotics of the modified Bessel functions on the Stokes lines}
\author{\sc R. B.\ Paris \\
{\em Division of Computing and Mathematics}, \\
{\em University of Abertay Dundee, Dundee DD1 1HG, UK}
}
\begin{document}
\def\f#1#2{\mbox{${\textstyle \frac{#1}{#2}}$}}
\def\dfrac#1#2{\displaystyle{\frac{#1}{#2}}}
\def\boldal{\mbox{\boldmath $\alpha$}}
{\newcommand{\Sgoth}{S\;\!\!\!\!\!/}
\newcommand{\bee}{\begin{equation}}
\newcommand{\ee}{\end{equation}}
\newcommand{\lam}{\lambda}
\newcommand{\ka}{\kappa}
\newcommand{\al}{\alpha}
\newcommand{\fr}{\frac{1}{2}}
\newcommand{\fs}{\f{1}{2}}
\newcommand{\g}{\Gamma}
\newcommand{\br}{\biggr}
\newcommand{\bl}{\biggl}
\newcommand{\ra}{\rightarrow}
\newcommand{\mbint}{\frac{1}{2\pi i}\int_{c-\infty i}^{c+\infty i}}
\newcommand{\mbcint}{\frac{1}{2\pi i}\int_C}
\newcommand{\mboint}{\frac{1}{2\pi i}\int_{-\infty i}^{\infty i}}
\newcommand{\gtwid}{\raisebox{-.8ex}{\mbox{$\stackrel{\textstyle >}{\sim}$}}}
\newcommand{\ltwid}{\raisebox{-.8ex}{\mbox{$\stackrel{\textstyle <}{\sim}$}}}
\renewcommand{\topfraction}{0.9}
\renewcommand{\bottomfraction}{0.9}
\renewcommand{\textfraction}{0.05}
\newcommand{\mcol}{\multicolumn}
\date{}
\maketitle
\pagestyle{myheadings}
\markboth{\hfill \sc R. B.\ Paris  \hfill}
{\hfill \sc  Asymptotics of the modified Bessel functions\hfill}
\begin{abstract}
We employ the exponentially improved asymptotic expansions of the confluent hypergeometric functions on the Stokes lines discussed by the author [Appl. Math. Sci. {\bf 7} (2013) 6601--6609] to give the analogous expansions of the modified Bessel functions $I_\nu(z)$ and $K_\nu(z)$ for large $z$ and finite $\nu$ on $\arg\,z=\pm\pi$ (and, in the case of $I_\nu(z)$, also on $\arg\,z=0$). Numerical results are presented to illustrate the accuracy of these expansions.

\vspace{0.4cm}

\vspace{0.3cm}

\end{abstract}

\vspace{0.3cm}

\noindent $\,$\hrulefill $\,$

\vspace{0.2cm}

\begin{center}
{\bf 1. \  Introduction}
\end{center}
\setcounter{section}{1}
\setcounter{equation}{0}
\renewcommand{\theequation}{\arabic{section}.\arabic{equation}}
The interest in exponentially precise asymptotics during the past three decades has shown that retention of exponentially small terms, previously neglected in asymptotics, can be essential for a high-precision description. 
For a discussion of recent developments in this area see \cite[Section 2.11]{DLMF} and \cite[Chapter 6]{PK}.
An example illustrating the advantage of retaining exponentially small terms in the asymptotic expansion of a certain integral is given in \cite[p.~66]{DLMF}. Although such terms are negligible in the Poincar\'e sense, their inclusion can significantly improve the numerical accuracy.

The modified Bessel function of the first kind $I_\nu(z)$ is defined by
\bee\label{e10}
I_\nu(z)=(\fs z)^\nu\sum_{k=0}^\infty\frac{(\fs z)^{2k}}{k! \g(k+\nu+1)}.
\ee
When $\nu$ takes on half-integer values, $I_\nu(z)$ can be expressed in terms of the hyperbolic functions; we exclude this special case from our asymptotic considerations. 
The behaviour of $I_\nu(z)$ for large $z$ and fixed $\nu$ is exponentially large throughout the sector $|\arg\,z|\leq\pi$, except on the imaginary axis $z=\pm ix$, $x>0$, where it is oscillatory (and is equal to $e^{\pm\pi\nu i/2}J_\nu(x)$, where $J_\nu(z)$ is the Bessel function of the first kind). The well-known asymptotic expansion of $I_\nu(z)$ for $|z|\to\infty$ and finite $\nu$ is given by \cite[p.~203]{WBF}, \cite[(10.40.5)]{DLMF} 
\bee\label{e11}
I_\nu(z)\sim\frac{e^z}{\sqrt{2\pi z}}\sum_{k=0}^\infty \frac{(-)^k a_k(\nu)}{z^k}+\frac{ie^{-z+\pi\nu i}}{\sqrt{2\pi z}}\sum_{k=0}^\infty \frac{a_k(\nu)}{z^k}
\ee
valid in the sector $-\fs\pi+\epsilon\leq\arg\,z\leq\f{3}{2}\pi-\epsilon$, where throughout $\epsilon$ denotes an arbitrary small constant. An analogous expansion with the sign of $i$ reversed holds in the conjugate sector $-\f{3}{2}\pi+\epsilon\leq\arg\,z\leq\fs\pi-\epsilon$. The coefficients $a_k(\nu)$ are defined by\footnote{This representation of the coefficients is equivalent to the familiar form $a_k(\nu)=(4\nu^2-1^2)(4\nu^2-3^2)\ldots(4\nu^2-(2k-1)^2)/(2^{3k} k!)$ for $k\geq 1$.}
\bee\label{e11a}
a_k(\nu)=\frac{(-)^k}{2^k k!}\,(\fs+\nu)_k (\fs-\nu)_k\qquad (k\geq 0),
\ee
where $(a)_k=\g(a+k)/\g(a)=a(a+1)\ldots (a+k-1)$ is the Pochhammer symbol. 

The first series in (\ref{e11}) is dominant as $|z|\to\infty$ in the right-half plane, with the second series being subdominant and maximally subdominant on $\arg\,z=0$. This situation is reversed in the left-half plane; the second series becomes dominant and assumes maximum
dominance over the first series on the Stokes lines $\arg\,z=\pm\pi$. 
On these rays the subdominant first series undergoes a Stokes phenomenon and begins to switch off (in the sense of increasing $|\arg\,z|$). An analogous process occurs on the other Stokes line $\arg\,z=0$, where the
subdominant second series in (\ref{e11}) undergoes a similar transition.
A correct interpretation of the expansion (\ref{e11}) would require the dominant series to be optimally truncated at, or near, its least term in magnitude (corresponding to the truncation index $k\simeq 2|z|$) in order to produce a remainder that is comparable to the subdominant contribution.

However, these Stokes transitions are not fully accounted for in the expansion (\ref{e11}). To see this, we put $z=xe^{\pi i}$, $x>0$, so that we have from (\ref{e11})
\bee\label{e12}
e^{-\pi\nu i} I_\nu(xe^{\pi i})\sim \frac{e^x}{\sqrt{2\pi x}}\sum_{k=0}^\infty \frac{(-)^k a_k(\nu)}{x^k}-\frac{ie^{-x-\pi\nu i}}{\sqrt{2\pi x}}\sum_{k=0}^\infty \frac{a_k(\nu)}{x^k}
\ee
as $x\to+\infty$. For real $\nu$, the dominant contribution in (\ref{e12}) is real, as it must be since from (\ref{e10}) we have 
\bee\label{e13}
I_\nu(xe^{\pm\pi i})=e^{\pm\pi\nu i} I_\nu(x)
\ee
with the consequence that $e^{-\pi\nu i}I_\nu(xe^{\pi i})$ is real in this case. However, the expansion (\ref{e12}) predicts a complex-valued exponentially small contribution (when $\nu$ is not a half-integer). The same reasoning applies to $I_\nu(x)$ on the Stokes line $x\in [0,\infty)$ since, from (\ref{e13}), its asymptotic expansion for $x\to+\infty$ is also given by the right-hand side of (\ref{e12}).

The modified Bessel function of the second kind $K_\nu(z)$ has the expansions
\begin{eqnarray*}
K_\nu(z)\sim\left\{\begin{array}{ll} \displaystyle{\sqrt{\frac{\pi}{2z}}\,e^{-z}\sum_{k=0}^\infty \frac{a_k(\nu)}{z^k}} & (|\arg\,z|\leq\f{3}{2}\pi\!-\!\epsilon) \\
\\
\displaystyle{\sqrt{\frac{\pi}{2z}}\,e^{-z}\sum_{k=0}^\infty \frac{a_k(\nu)}{z^k}+2i\cos \pi\nu\,\sqrt{\frac{\pi}{2z}}\,e^{z}\sum_{k=0}^\infty \frac{(-)^ka_k(\nu)}{z^k}}&\\
&\!\!\!\!\!\!\!\!(\fs\pi\!+\!\epsilon\leq\arg\,z\leq \f{5}{2}\pi\!-\!\epsilon).\end{array}\right.
\end{eqnarray*}
Consequently, on the Stokes line $\arg\,z=\pi$, we obtain with $z=xe^{\pi i}$, $x>0$
\bee\label{e14}
K_\nu(xe^{\pi i})\sim -i\sqrt{\frac{\pi}{2x}}\,e^x\sum_{k=0}^\infty \frac{(-)^ka_k(\nu)}{x^k}+2\cos \pi\nu\,\sqrt{\frac{\pi}{2x}}\,e^{-x}\sum_{k=0}^\infty \frac{a_k(\nu)}{x^k}
\ee
as $x\to+\infty$. A conjugate expansion holds on $\arg\,z=-\pi$. However, the expansion (\ref{e14}) will not yield the correct form for the exponentially small contribution for large $x$, since on $\arg\,z=\pi$ the appearance of the subdominant second series in (\ref{e14}) is only halfway through its Stokes transition.

In this note we employ the exponentially improved asymptotic expansions of the confluent hypergeometric functions on the Stokes lines discussed in \cite{PCHF} to deal with the analogous expansions of the modified Bessel functions.
We obtain the exponentially small contribution associated with $I_\nu(z)$ on the Stokes lines $\arg\,z=0$ and $\arg\,z=\pi$, and that associated with $K_\nu(z)$ on the Stokes line $\arg\,z=\pi$, that take into account the incomplete nature of the Stokes transitions on these rays. 

\vspace{0.6cm}

\begin{center}
{\bf 2. \ A summary of the asymptotic expansions of the Kummer functions}
\end{center}
\setcounter{section}{2}
\setcounter{equation}{0}
\renewcommand{\theequation}{\arabic{section}.\arabic{equation}}
We present in this section a summary of the exponentially improved expansion of the Kummer functions ${}_1F_1(a;b;-x)$ and $U(a,b,xe^{\pm\pi i})$ for $x\to+\infty$ with fixed parameters $a$ and $b$.

The function $U(a,b,z)$ has, in general, a branch point at $z=0$
with the $z$-plane cut along $(-\infty,0]$ and possesses the large-$|z|$ expansion
\[U(a, b, z) \sim z^{-a} \sum_{k=0}^\infty\frac{(-)^k(a)_k (1+a-b)_k}{k!\,z^k}\qquad (|\arg\,z|\leq\f{3}{2}\pi-\epsilon),\]
which is algebraic in character in the stated sector. However, on the Stokes lines $\arg\,z=\pm\pi$, an exponentially small expansion switches on, so that in the sectors 
$\pi\leq|\arg\,z|\leq\f{3}{2}\pi-\epsilon$ we have a compound expansion with a subdominant exponential contribution. This latter contribution becomes dominant beyond $|\arg\,z|=\f{3}{2}\pi$. 

In the following we let $M$, $N$, $m$ denote positive integers  and define the parameters 
\bee\label{e20}
\vartheta:=a-b,\qquad \mu\equiv\mu(m):=2a-b+m.
\ee
The exponentially improved expansion of $U(a,b,z)$ is given by \cite[p.~329]{DLMF}
\[U(a, b, z)=z^{-a}\sum_{k=0}^{m-1}\frac{(-)^k(a)_k (1+a-b)_k}{k!\,z^k}\hspace{6cm}\]
\bee\label{e21}
\hspace{4cm}+\frac{2\pi i e^{-\pi i(a+\vartheta)}}{\g(a) \g(1+a-b)}\,z^\vartheta e^{z}\bl\{\sum_{j=0}^{M-1}A_j z^{-j}\,T_{\nu-j}(z)+R_{M,m}(z)\br\},
\ee
where 
\bee\label{e21a}
A_j=\frac{(1-a)_j(1-b)_j}{j!}\qquad(j\geq 0)
\ee
and $T_\mu(z)$ denotes the so-called {\it terminant function} defined as a multiple of the 
incomplete gamma function $\g(a,z)$ by
\[T_\mu(z):=\frac{e^{\pi i\mu}\g(\mu)}{2\pi i}\,\g(1-\mu,z).\]
In (\ref{e21}),
$m$ is an arbitrary positive integer but will be chosen to be the optimal truncation index $m_o$ of the algebraic expansion corresponding to truncation at, or near, the least term in magnitude. This is easily verified to be 
\bee\label{e2o}
m_o\simeq |z|-\Re (2a-b),
\ee
so that $m_o\ra\infty$ as $|z|\ra+\infty$. When $m-|z|$ is bounded, the remainder term
in (\ref{e21}) satisfies $R_{M,m}(z)=O(e^{-|z|-z}z^{-M})$ as $|z|\ra\infty$ in the sector
$|\arg\,z|\leq\pi$. 

Since the truncation index $m$ is chosen to be optimal, the index $\mu$ appearing in (\ref{e20}) satisfies  $\mu\sim |z|$ as $|z|\ra+\infty$. The asymptotic expansion of $T_\mu(z)$ for large $|\mu|$ and $|z|$, when $|\mu|\sim |z|$, has been discussed in detail by Olver in \cite{O}. By expressing $T_\mu(z)$ in terms of a Laplace integral, which is associated with a saddle point and a simple pole becoming coincident on $\arg\,z=\pi$, Olver \cite[\S 5]{O} established that, for $z=xe^{\pi i}$, $x>0$,
\bee\label{e26}
T_{\mu-j}(xe^{\pi i})=\frac{1}{2}-\frac{i}{\sqrt{2\pi x}}\bl\{\sum_{k=0}^{N-1} (\fs)_{k} g_{_{2k}}(j) (\fs x)^{-k}+O(x^{-N})\br\}\qquad (x\ra+\infty),
\ee
where the coefficients $g_{_k}(j)$ result from the expansion 
\[\frac{\tau^{\gamma_j-1}}{1-\tau}\,\frac{d\tau}{dw}=-\frac{1}{w}+\sum_{k=0}^\infty g_{_k}(j)w^k,\qquad \fs w^2=\tau-\log\,\tau-1.\]

The branch of $w(\tau)$ is chosen such that $w\sim \tau-1$ as $\tau\ra 1$ and the parameter $\gamma_j$ is specified by
\bee\label{e26a}
\gamma_j=\mu-x-j\qquad(0\leq j\leq N-1)
\ee
with $|\gamma_j|$ bounded. Upon reversion of the $w$-$\tau$ mapping to yield
\[\tau=1+w+\f{1}{3}w^2+\f{1}{36}w^3-\f{1}{270}w^4+\f{1}{4320}w^5+ \cdots\ ,\]
it is found with the help of {\it Mathematica} that the first five even-order coefficients $g_{_{2k}}(j)\equiv 6^{-2k} {\hat g}_{_{2k}}(j)$ are\footnote{There was a misprint in the first term in ${\hat g}_6(j)$ in \cite{PCHF}, which appeared as $-3226$ instead of $-3626$. This was pointed out by T. Pudlik \cite{TP}. The correct value was used in the numerical calculations described in \cite{PCHF}.}
\begin{eqnarray*}
{\hat g}_0(j)\!\!&=&\!\!\f{2}{3}-\gamma_j,\qquad {\hat g}_2(j)=\f{1}{15}(46-225\gamma_j+270\gamma_j^2-90\gamma_j^3), \\
{\hat g}_4(j)\!\!&=&\!\!\f{1}{70}(230-3969\gamma_j+11340\gamma_j^2-11760\gamma_j^3+5040\gamma_j^4
-756\gamma_j^5),\\
{\hat g}_6(j)\!\!&=&\!\!\f{1}{350}(-3626-17781\gamma_j+183330\gamma_j^2-397530\gamma_j^3+370440\gamma_j^4
-170100\gamma_j^5\\
&&\hspace{7cm}+37800\gamma_j^6-3240\gamma_j^7),\\
{\hat g}_8(j)\!\!&=&\!\!\f{1}{231000}(-4032746+43924815\gamma_j+88280280\gamma_j^2-743046480\gamma_j^3\\
&&+1353607200\gamma_j^4-1160830440\gamma_j^5+541870560\gamma_j^6
-141134400\gamma_j^7\\
&&\hspace{6cm}+19245600\gamma_j^8-1069200\gamma_j^9).
\end{eqnarray*}

Substitution of (\ref{e26}) into (\ref{e21}) (where for convenience we put $M=N$) 
and introduction of the coefficients $B_j$ defined by
\bee\label{e28}
B_j=\sum_{k=0}^j (-2)^{k} (\fs)_k\,A_{j-k}\,g_{_{2k}}(j-k),
\ee
then yields the expansion: 
\newtheorem{theorem}{Theorem}
\begin{theorem}$\!\!\!. \ \cite[(3.2)]{PCHF}$
We have the expansion 
\[U(a, b, xe^{\pm\pi i})-(xe^{\pm\pi i})^{-a}\sum_{k=0}^{m_o-1}\frac{(a)_k (1+a-b)_k}{k!\,x^k}\hspace{8cm}\]
\bee\label{e29}
=\pm\frac{2\pi i\, e^{\mp\pi ia}\,x^{a-b} e^{-x}}{\g(a) \g(1+a-b)}\bl\{\frac{1}{2}\sum_{j=0}^{M-1}(-)^jA_j x^{-j}\mp\frac{i}{\sqrt{2\pi x}}\sum_{j=0}^{M-1}(-)^jB_j x^{-j}+O(x^{-M})\br\}
\ee
as $x\ra+\infty$,  provided $a,\,1+a-b\neq 0, -1, -2, \ldots\,$. The integer $m_o$ is the optimal truncation index of the algebraic expansion satisfying $m_o\sim x$, $M$ is a positive integer and the coefficients $A_j$ and $B_j$ are defined in (\ref{e21a}) and (\ref{e28}), respectively.
\end{theorem}

From the relation connecting the first Kummer function ${}_1F_1(a; b; -x)$ to $U(a,b,xe^{\pm\pi i})$ \cite[p.~323]{DLMF} we obtain the expansion:
\begin{theorem}$\!\!\!.\ \cite[(2.11)]{PCHF}$
When $\vartheta=a-b$ is non-integer we have the expansion 
\[\frac{\g(a)}{\g(b)}\,{}_1F_1(a; b; -x)-\frac{x^{-a}\g(a)}{\g(b-a)} \sum_{k=0}^{m_o-1} \frac{(a)_k (1+a-b)_k}{k!\,x^k}\hspace{6cm}\]
\bee\label{e210}
\hspace{2cm}=x^{a-b} e^{-x}\bl\{ \cos \pi\vartheta \sum_{j=0}^{M-1}(-)^jA_jx^{-j}-\frac{2\sin \pi\vartheta}{\sqrt{2\pi x}}\sum_{j=0}^{M-1}(-)^j B_jx^{-j}+O(x^{-M})\br\}
\ee
as $x\ra+\infty$, where $m_o$ is the optimal truncation index of the algebraic expansion satisfying $m_o\sim x$, $M$ is a positive integer and the coefficients $A_j$ and $B_j$ are defined in (\ref{e21a}) and (\ref{e28}), respectively.
\end{theorem}
When $\vartheta=n$, $n=0, 1, 2, \ldots\ $, the algebraic expansion in (\ref{e210}) vanishes and the coefficients $A_j$ vanish for $j>n$. In this case the function ${}_1F_1(a;b;z)$ is a polynomial in $z$ of degree $n$. When $\vartheta=-n$, $n=1, 2, \ldots\ $, the algebraic expansions in (\ref{e29}) and (\ref{e210}) consist of $n$ terms (and so cannot be optimally truncated); see \cite{PCHF} for details. In both cases the second exponentially small series in (\ref{e210}) vanishes.
\vspace{0.6cm}

\begin{center}
{\bf 3. \ The expansions on the Stokes lines}
\end{center}
\setcounter{section}{3}
\setcounter{equation}{0}
\renewcommand{\theequation}{\arabic{section}.\arabic{equation}}
From the identity expressing $I_\nu(z)$ in terms of the confluent hypergeometric function \cite[(10.39.5)]{DLMF}
we have, with $z=xe^{\pi i}$, $x>0$
\[e^{-\pi\nu i}I_\nu(xe^{\pi i})=\frac{(\fs x)^\nu e^x}{\g(\nu+1)}\,{}_1F_1(\nu+\fs; 2\nu+1; -2x).\]
With the parameter $\vartheta=-\nu-\fs$, we obtain from Theorem 2 the expansion
\[e^{-\pi\nu i}I_\nu(xe^{\pi i})=\frac{(\fs x)^\nu e^x}{\g(\nu+1)}\,\frac{\g(2\nu+1)}{\g(\nu+1)}\left\{(2x)^{-\nu-\fr}\sum_{k=0}^{m_o-1}\frac{(\fs+\nu)_k (\fs-\nu)_k}{k! (2x)^k}\right.\]
\[\left.+(2x)^{-\nu-\fr}e^{-2x}\bl\{\cos \pi\vartheta \sum_{j=0}^{M-1}\frac{(-)^j A_j}{(2x)^j}-\frac{2\sin \pi\vartheta}{\sqrt{4\pi x}}\sum_{j=0}^{M-1}\frac{(-)^j B_j}{(2x)^j}+O(x^{-M})\br\}\right\}.\]
Use of the duplication formula for the gamma function
\[\g(2z)=2^{2z-1}\pi^{-\fr} \g(z) \g(z+\fs)\]
and the fact that the coefficients $a_j(\nu)=(-2)^{-j} A_j$, where $A_j$ are defined in (\ref{e21a}) with $a=\nu+\fs$, $b=2\nu+1$, then produces
\begin{theorem}$\!\!\!.$ Let $\vartheta=-\nu-\fs$ and $M$ be a positive integer. Then we have the expansion
\[
e^{-\pi\nu i}I_\nu(xe^{\pi i})=\frac{e^x}{\sqrt{2\pi x}} \sum_{k=0}^{m_o-1}\frac{(-)^ka_k(\nu)}{x^k}\hspace{7cm}\]
\bee\label{e31}
\hspace{3cm}+\frac{e^{-x}}{\sqrt{2\pi x}}\bl\{\cos \pi\vartheta \sum_{k=0}^{M-1}\frac{a_k(\nu)}{x^k}-\frac{\sin \pi\vartheta}{\sqrt{\pi x}}\sum_{j=0}^{M-1}\frac{(-)^j B_j}{(2x)^j}+O(x^{-M})\br\}
\ee
as $x\to+\infty$, where $m_o$ is the optimal truncation index for the dominant series  satisfying $m_o\simeq 2x$.
The coefficients $a_k(\nu)$ and $B_j$ are defined in (\ref{e11a}) and (\ref{e28}), respectively.
\end{theorem}
We remark that the right-hand side of (\ref{e31}) also gives the expansion of $e^{\pi\nu i}I_\nu(xe^{-\pi i})$ as well as that of $I_\nu(x)$ as $x\to+\infty$ since,  
by (\ref{e13}), $I_\nu(x)=e^{\mp\pi\nu i}I_\nu(xe^{\pm\pi i})$. 

The modified Bessel function $K_\nu(z)$ is given by
\[K_\nu(z)=(2z)^\nu \sqrt{\pi} e^{-z}\,U(\nu+\fs, 2\nu+1, 2z).\]
From Theorem 1 we then find\footnote{This result can also be obtained from (10.40.2) and (10.40.13) in \cite{DLMF} combined with the expansion (\ref{e26}).} with $z=xe^{\pm\pi i}$, $x>0$ that
\[K_\nu(xe^{\pm\pi i})=(2xe^{\pm\pi i})^\nu \sqrt{\pi} e^x\left\{(2xe^{\pm\pi i})^{-\nu-\fr}\sum_{k=0}^{m_o-1}\frac{(\fs+\nu)_k(\fs-\nu)_k}{k! (2x)^k}\right.\]
\[\left.\pm \frac{2\pi i\,(2xe^{\pm\pi i})^{-\nu-\fr}e^{-2x}}{\g(\fs+\nu) \g(\fs-\nu)}\bl\{\frac{1}{2}\sum_{j=0}^{M-1}\frac{(-)^jA_j}{(2x)^j}\mp\frac{i}{2\sqrt{\pi x}}\sum_{j=0}^{M-1}\frac{(-)^jB_j}{(2x)^j}+O(x^{-M})\br\}\right\}.
\]
After some straightforward rearrangement, this produces
\begin{theorem}$\!\!\!.$ Let $M$ be a positive integer. Then we have the expansions
\[K_\nu(xe^{\pm\pi i})=\mp i\sqrt{\frac{\pi}{2x}}\,e^x\sum_{k=0}^{m_o-1} \frac{(-)^k a_k(\nu)}{x^k}\hspace{7cm}\]
\bee\label{e32}
+2\cos \pi\nu\,\sqrt{\frac{\pi}{2x}} e^{-x}\bl\{\frac{1}{2}\sum_{k=0}^{M-1}\frac{a_k(\nu)}{x^k}\mp\frac{i}{2\sqrt{\pi x}} \sum_{j=0}^{M-1} \frac{(-)^jB_j}{(2x)^j}+O(x^{-M})\br\}
\ee
as $x\to+\infty$, where $m_o$ is the optimal truncation index for the dominant series  satisfying $m_o\simeq 2x$.
The coefficients $a_k(\nu)$ and $B_j$ are defined in (\ref{e11a}) and (\ref{e28}), respectively.
\end{theorem}
Comparison of (\ref{e32}) with (\ref{e14}) reveals the by-now familiar fact that the value of the Stokes multiplier of $K_\nu(z)$  on $\arg\,z=\pm\pi$ (given by the expression in curly braces in (\ref{e32})) is $\fs$ to leading order.

In numerical calculations, we set the optimal truncation index $m_o=2x+\alpha$, where $|\alpha|\leq\fs$; when $2x$ is an integer then $\alpha=0$. Then from (\ref{e20}) and (\ref{e26a}) we have
\[\gamma_j=m_o-2x-j=\alpha-j.\]
The coefficients $B_j$ can be computed from (\ref{e28}) and are real when $\nu$ is real. We subtract off the dominant, exponentially large series in (\ref{e31}) by defining
\[F_\nu(x):=e^{-\pi\nu i}I_\nu(xe^{\pi i})-\frac{e^x}{\sqrt{2\pi x}} \sum_{k=0}^{m_o-1}\frac{(-)^ka_k(\nu)}{x^k}~.\]
The exponentially small expansion is given by
\[S_I(M;x):=\frac{e^{-x}}{\sqrt{2\pi x}}\bl\{\cos \pi\vartheta \sum_{k=0}^{M-1}\frac{a_k(\nu)}{x^k}-\frac{\sin \pi\vartheta}{\sqrt{\pi x}} \sum_{j=0}^{M-1}\frac{(-)^j B_j}{(2x)^j}\br\},\]
which is seen to be real when $\nu$ is real.
In Table 1 we show\footnote{In Tables 1 and 2 we write the values as $x(y)$ instead of $x\times 10^y$.} the values of $S_I(M;x)$ as a function of the truncation index $M$ for different values of $x$ and compare these with the value of $F_\nu(x)$.
\begin{table}[th]
\caption{\footnotesize{The values of $F_\nu(x)$ and $S_I(M;x)$ for $e^{-\pi\nu i}I_\nu(xe^{\pi i})$ for $x=25$ and different truncation index $M$ when $\nu=1/4$. }}
\begin{center}
\begin{tabular}{|c|l|l|l|}
\hline
&&&\\[-0.3cm]
\mcol{1}{|c|}{$M$} & \mcol{1}{c|}{$S_I(M;10)$} & \mcol{1}{c|}{$S_I(M;15.4)$} & \mcol{1}{c|}{$S_I(M;20)$}\\
[.1cm]\hline
&&&\\[-0.3cm]
1 & $-3.568247262(-06)$ & $-1.163196884(-08)$ & $-1.190793339(-10)$\\
2 & $-3.538491386(-06)$ & $-1.151737379(-08)$ & $-1.185631040(-10)$\\
3 & $-3.539961940(-06)$ & $-1.151900437(-08)$ & $-1.185763303(-10)$\\
4 & $-3.539827969(-06)$ & $-1.151884096(-08)$ & $-1.185756985(-10)$\\
5 & $-3.539846361(-06)$ & $-1.151885440(-08)$ & $-1.185757438(-10)$\\
6 & $-3.539842998(-06)$ & $-1.151885272(-08)$ & $-1.185757394(-10)$\\
7 & $-3.539843764(-06)$ & $-1.151885298(-08)$ & $-1.185757400(-10)$\\
[.1cm]\hline
&&&\\[-0.3cm]
$F_\nu(x)$ & $-3.539843604(-06)$ & $-1.151885294(-08)$ & $-1.185757399(-10)$\\
[.1cm]\hline
\end{tabular}
\end{center}
\end{table}

Similarly, we define
\[G_\nu(x):=K_\nu(xe^{\pi i})+i\sqrt{\frac{\pi}{2x}}\,e^x \sum_{k=0}^{m_o-1}\frac{(-)^k a_k(\nu)}{x^k}\]
and
\[S_K(M;x):=2\cos \pi\nu\,\sqrt{\frac{\pi}{2x}} e^{-x}\bl\{\frac{1}{2}\sum_{k=0}^{M-1}\frac{a_k(\nu)}{x^k}-\frac{i}{2\sqrt{\pi x}} \sum_{j=0}^{M-1} \frac{(-)^jB_j}{(2x)^j}\br\}.\]
Table 2 shows an example of the values of $S_K(M;x)$ for different truncation index $M$ compared with the value of $G_\nu(x)$. It can be seen in both cases that the computed values of $F_\nu(x)$ and $G_\nu(x)$ agree well their corresponding exponentially small expansions.
\begin{table}[h]
\caption{\footnotesize{The values of $G_\nu(x)$ and $S_K(M;x)$ for $K_\nu(xe^{\pi i})$ for different truncation index $M$ when $x=25$ and $\nu=1/4$. }}
\begin{center}
\begin{tabular}{|c|l|}
\hline
&\\[-0.3cm]
\mcol{1}{|c|}{$M$} & \mcol{1}{c|}{$S_K(M;25)$} \\
[.1cm]\hline
&\\[-0.3cm]
1 & $2.461573958(-12) - 1.851725849i(-13)$\\
2 & $2.452343056(-12) - 1.839098107i(-13)$\\
3 & $2.452544982(-12) - 1.839470730i(-13)$\\
4 & $2.452536653(-12) - 1.839451010i(-13)$\\
5 & $2.452537160(-12) - 1.839452410i(-13)$\\
6 & $2.452537119(-12) - 1.839452283i(-13)$\\
7 & $2.452537123(-12) - 1.839452297i(-13)$\\
[.1cm]\hline
&\\[-0.3cm]
$G_\nu(x)$ & $2.452537123(-12)-1.839452296i(-13)$ \\
[.1cm]\hline
\end{tabular}
\end{center}
\end{table}


\vspace{0.6cm}

\begin{thebibliography}{99}
\footnotesize{

\bibitem{O}
F.W.J. Olver, Uniform, exponentially improved asymptotic expansions for the generalized exponential integral, SIAM J. Math. Anal. {\bf 22} (1991) 1460--1474.


\bibitem{DLMF}
F.W.J. Olver, D.W. Lozier, R.F. Boisvert and C.W. Clark (eds.),    
{\it NIST Handbook of Mathematical Functions}, Cambridge University Press, Cambridge, 2010.

\bibitem{PK} 
R.B. Paris and D. Kaminski,  {\it Asymptotics and Mellin-Barnes Integrals}, 
Cambridge University Press, Cambridge, 2001.


\bibitem{PCHF}
R.B. Paris, Exponentially small expansions of the confluent hypergeometric functions, Appl. Math. Sci. {\bf 7} (2013) 6601--6609.

\bibitem{TP}
T. Pudlik, Private communication (2016).

\bibitem{WBF}
G.N. Watson, {\it Theory of Bessel Functions}, Cambridge University Press, Cambridge, 1952. 


}
\end{thebibliography}
\end{document}